\documentclass[11pt]{amsart}
\input epsf
\usepackage{latexsym}
\usepackage{amssymb,amsmath,amsthm,amscd,amssymb,enumerate, verbatim}
\usepackage[all]{xy}

\numberwithin{equation}{section}

\newtheorem{thm}[equation]{Theorem}

\newtheorem{Example}[equation]{Example}

\newtheorem{remark}[equation]{Remark}
\newenvironment{rmk}{\begin{remark}\rm}{\end{remark}}
\def\co{\colon\thinspace}

\newcommand{\Int}{\mbox{Int}}

\def\a{\alpha}

\def\b{\beta}

\def\d{\partial}

\def\S1{\bf S^1}

\begin{document}

\title[Open aspherical manifolds not covered by Euclidean space]
{\bf Open aspherical manifolds not covered by the Euclidean space}
\thanks{
{\it\! Keywords:\rm\ aspherical manifold, homology sphere, 
hyperplane unknotting}.\rm\
The author is grateful for NSF support (DMS-1105045)}\rm


\author{Igor Belegradek}

\address{Igor Belegradek\\School of Mathematics\\ Georgia Institute of
Technology\\ Atlanta, GA 30332-0160}\email{ib@math.gatech.edu}

\date{}

\begin{abstract} 
We show that any open aspherical manifold of dimension $n\ge 4$
is tangentially homotopy equivalent to an $n$-manifold
whose universal cover is not homeomorphic to $\mathbb R^n$.
\end{abstract}

\maketitle

\section{Introduction}

Davis famously constructed, for every $n\ge 4$, a
closed aspherical $n$-manifold whose universal cover is 
not homeomorphic to $\mathbb R^n$~\cite{Dav-annals}. 
We prove:

\begin{thm}
\label{thm: univ cover not Rn}
If the universal cover of 
an open smooth manifold $V$ is diffeomorphic to
$\mathbb R^n$ with $n\ge 4$, then the tangential homotopy type of $V$ contains a continuum of open 
smooth $n$-manifolds whose universal covers are not homeomorphic to $\mathbb R^n$. 
\end{thm}

Our proof also works in {\scshape pl} and {\scshape top} categories.
Recall that a homotopy equivalence $f\co V\to W$ of $n$-manifolds
is {\it tangential\,} if $f^\#TW$  and $TV$ are stably isomorphic,
where $TX$ denotes the tangent bundle if
$X$ is a smooth manifold and the tangent microbundle in the 
{\scshape pl} or {\scshape top} cases. 
By~\cite[Theorem 2.3]{Sie-collar} this is equivalent to
requiring that $f\times\mathbf{id}(\mathbb R^s)$ is
homotopic to a {\scshape cat}-homeomorphism for some $s$
(where {\scshape cat} equals {\scshape diff}, {\scshape pl}
or {\scshape top}).

In the simply-connected case Theorem~\ref{thm: univ cover not Rn} 
is due to Curtis-Kwun~\cite{CurKwu65} 
for $n\ge 5$ and
to Glasner~\cite{Gla67} for $n=4$.
The proof combines three ingredients: 

(1) A result of Curtis-Kwun~\cite{CurKwu65} that for
a boundary connected sum $S$ of a countable family 
of compact $n$-manifolds, the homeomorphism type
of $\Int(S)$ determines the isomorphism class of 
$\pi_1(\d S)$.

(2)
A recent result
of Calcut-King-Siebenmann~\cite{CKS} that
any countable collection
of {\scshape cat} properly embedded $\mathbb R^{n-1}$'s in $\mathbb R^n$ is 
{\scshape cat} unknotted, which generalizes classical 
results of Cantrell and Stallings.

(3) The existence of infinitely many smooth compact contractible
$n$-manifolds whose boundary homology spheres
have freely indecomposable fundamental groups.
(Such examples can be found 
in~\cite{CurKwu65, Gla67}, and more examples
are now known, see e.g. Casson-Harer~\cite{CasHar81} 
for aspherical homology $3$-spheres that bound
smooth contractible $4$-manifolds, while 
Kervaire~\cite{Ker69} showed that
the fundamental group of any homology $3$-sphere appears
as the fundamental group of the boundary of a smooth
contractible $n$-manifolds for any $n\ge 5$).

\begin{proof}
Since $V$ is open, it contains a {\scshape cat} properly embedded ray
whose {\scshape cat} regular neighborhood is an embedded closed 
halfspace, see e.g.~\cite[Section 3]{CKS}. 
Hence $V$ is {\scshape cat} isomorphic
to the interior of a noncompact
manifold $N$ whose boundary is an open disk.
By the strong version of
the Cantrell-Stalling hyperplane unknotting 
theorem proved in~\cite[Corollary 9.3]{CKS}, 
the universal cover of $N$
can be compactified to $D^n$, the $n$-disk, 
in which the preimage of $\d N$ becomes
the union of a countable collection of 
round open disks with pairwise disjoint closures. 

Let $\{C_i\}_{i\in\mathbb N}$ be an
infinite sequence of compact contractible $n$-manifolds,
such that $\pi_1(\d C_i)$ are pairwise non-isomorphic
and freely indecomposable.
Given a subset $\a\subseteq\mathbb N$, let $C_\a$ be
a boundary connected sums of $C_i$'s
with indices in $\a$. 
(For our purposes the
choices involved in defining 
boundary connected sums will always be irrelevant).
Fix a closed $(n-1)$-disk $\Delta\subset \d C_\a$,
and let $N_\a$ be a boundary connected sum
of $N$ and $C_\a$ obtained by identifying $\Delta$
with a closed disk in $\d N$. 
A deformation retraction $C_\a\to\Delta$ extends to
a deformation retraction of $N_\a\to N$, so 
$V_\a:=\Int(N_\a)$ is 
tangentially homotopy equivalent to $V$.

If $Q_\a$ denotes a boundary
connected sum of countably many copies of $C_\a$, then
the interior of the universal cover of $N_\a$
is  homeomorphic to the interior of a boundary connected sum 
of $D^n$ and $Q_\a$, which
is homeomorphic to $\Int(Q_\a)$.
By~\cite[Theorem 4.1]{CurKwu65}
if $Q_\a$, $Q_\b$ have homeomorphic interiors,
then $\d Q_\a$, $\d Q_\b$ have isomorphic fundamental groups.
Now $\pi_1(\d Q_\a)$ is a free product
in which each factor $\pi_1(\d C_{i_k})$, $i_k\in\a$ appears
countably many times.
Each $\pi_1(\d C_{i_k})$ is freely indecomposable,
so $\a=\b$ by Grushko's theorem.
Thus the universal covers of
$\Int(N_\a)$ lie in a continuum of homeomorphism types.
\end{proof}

\begin{rmk}
The proof of~\cite[Theorem 4.1]{CurKwu65}
is quite technical, which may be due to the fact that
the tools of Siebenmann's thesis
were not yet available at the time, so we
summarize it in a modern language: 
Given $\a=\{i_1, \dots, i_k, \dots\}$
it is easy to construct
a cofinal family $\{U_k\}_{k\ge 1}$ 
of neighborhoods of infinity in $Q_\a$
such that in the corresponding inverse 
sequence of fundamental groups,
the group $\{\pi_1(U_k)\}$ 
is the free product of 
$\pi_1(C_{i_1})\ast\dots\ast\pi_1(C_{i_k})$,
and the map $\pi_1(U_k)\longleftarrow \pi_1(U_{k+1})$
is a retraction onto the first $k$ factors,
and in particular, is surjective. Hence the inverse
sequence of groups is Mittag-Leffler, and therefore
its pro-equivalence
class depends only on the homeomorphism
type of $Q_\a$. Now if $Q_\a$, $Q_\b$ have 
homeomorphic interiors, then a simple diagram
chase in the commutative diagram from the definition
of pro-equivalence shows that each free factor
of $\pi_1(Q_\a)$ occurs as a free factor of $\pi_1(Q_\b)$,
so $\a=\b$.
\end{rmk}

\begin{rmk} 
Theorem~\ref{thm: univ cover not Rn} should hold for $n=3$,
but our proof fails. One could try
substituting the boundary connected sum of $N$ and $Q_\a$
by the end sum of $V$ with 
a suitable Whitehead manifold, but 
the multiple hyperplane unknotting is no longer true,
due to existence of an exotic 
$[0,1]\times \mathbb R^2$~\cite{ScoTuc}. 
This makes analyzing the fundamental group
of infinity more delicate. For the same reason
we do not attempt to prove Theorem~\ref{thm: univ cover not Rn}
for $V$ whose universal cover is not $\mathbb R^n$.
\end{rmk}

\begin{rmk}
By the Cartan-Hadamard theorem
any (complete Riemannian)
$n$-manifold of nonpositive curvature is covered by
$\mathbb R^n$. So given an open 
$n$-manifold of nonpositive curvature with $n\ge 4$,
Theorem~\ref{thm: univ cover not Rn} 
yields a continuum of $n$-manifolds 
in the same tangential homotopy type
that admit no metric of nonpositive curvature.
In~\cite{Bel-bus} the
author studied obstructions to nonpositive
curvature on open manifolds that go beyond the Cartan-Hadamard theorem.
\end{rmk}

\small
\bibliographystyle{amsalpha}

\def\cprime{$'$}
\providecommand{\bysame}{\leavevmode\hbox to3em{\hrulefill}\thinspace}
\providecommand{\MR}{\relax\ifhmode\unskip\space\fi MR }
\providecommand{\MRhref}[2]{%
  \href{http://www.ams.org/mathscinet-getitem?mr=#1}{#2}
}
\providecommand{\href}[2]{#2}

\end{document}